\documentclass[11pt,a4paper]{article}
\usepackage{amsmath}
\usepackage{amsfonts}
\usepackage{amssymb}
\usepackage{graphicx}
\usepackage{enumerate}

\newtheorem{thm}{Theorem}[section]
\newtheorem{prop}[thm]{Proposition}

\newtheorem{con}[thm]{Conjecture}
\newtheorem{lemma}[thm]{Lemma}

\newenvironment{proof} {\par \noindent \textbf{Proof: }}{\QED \par \bigskip \par}
\newcommand{\QED}{\hfill$\square$}
\newcommand{\rz}{\vspace{0.2cm}}

\begin{document}
\baselineskip 16pt

\phantom{.} \vskip 6cm

\begin{center}
{\LARGE \bf ON COMPARING ZAGREB INDICES \footnote{%
Supported by Research Grants 144007 and 144015G of Serbian Ministry
of Science and Environmental Protection, and Research Program
P1-0285 of Slovenian Agency for Research. } }

\bigskip
\bigskip
{\large \sc Aleksandar Ili\'c}

\medskip
{\em \normalsize Faculty of Sciences and Mathematics, University of Ni\v s} \\
{\em \normalsize Department of Mathematics and Informatics,  Vi\v
segradska 33, 18000 Ni\v s, Serbia}\\
{\normalsize e-mail: { \tt aleksandari@gmail.com }}

\bigskip
{\large \sc Dragan Stevanovi\'c}

\medskip
{\em \normalsize University of Primorska---FAMNIT, Glagolja\v ska 8, 6000 Koper, Slovenia,} \\
{\em \normalsize Mathematical Institute, Serbian Academy of Science
and Arts,}\\
{\em \normalsize Knez Mihajlova 36, 11000 Belgrade, Serbia }\\
{\normalsize e-mail: { \tt dragance106@yahoo.com}}

\bigskip\medskip
{\small (Received October 21, 2008)}
\bigskip
\end{center}

\begin{abstract}
Let $G=(V,E)$ be a simple graph with $n = |V|$ vertices and $m =
|E|$ edges. The first and second Zagreb indices are among the oldest
and the most famous topological indices, defined as $M_1 = \sum_{i
\in V} d_i^2$ and $M_2 = \sum_{(i, j) \in E} d_i d_j$, where $d_i$
denote the degree of vertex~$i$. Recently proposed conjecture $M_1 /
n \leqslant M_2 / m$ has been proven to hold for trees, unicyclic
graphs and chemical graphs, while counterexamples were found for
both connected and disconnected graphs. Our goal is twofold, both in
favor of a conjecture and against it. Firstly, we show that the
expressions $M_1/n$ and $M_2/m$ have the same lower and upper
bounds, which attain equality for and only for regular graphs. We
also establish sharp lower bound for variable first and second
Zagreb indices. Secondly, we show that for any fixed number
$k\geqslant 2$, there exists a connected graph with $k$~cycles for
which $M_1/n>M_2/m$ holds, effectively showing that the conjecture
cannot hold unless there exists some kind of limitation on the
number of cycles or the maximum vertex degree in a graph. In
particular, we show that the conjecture holds for subdivision
graphs.
\end{abstract}


\section{Introduction}

Let $G = (V, E)$ be a simple graph with $n = |V|$ vertices and $m =
|E|$ edges. The first Zagreb index $M_1$ and the second Zagreb index
$M_2$ of $G$ are defined as follows:
$$
M_1 = \sum_{i \in V} d_i^2 \quad \mbox { and } \quad M_2 = \sum_{(i,
j) \in E} d_i d_j,
$$
where $d_1, d_2, \ldots, d_n$ are vertex degrees, while $d_i d_j$
represents weight associated to the edge $(i,j)$. The Zagreb indices
were first introduced in~\cite{GuTr72} and the survey of properties
of $M_1$ and $M_2$ is given in~\cite{NiKMT03}. Note that in random
graphs with $n$ vertices and uniform edge probability~$p$, the order
of magnitude of $M_{1}$ is $O(n^{3}p^{2})$, while the order of
magnitude of $M_{2}$ is $O(n^{4}p^{3})$, implying that $M_{1}/n$ and
$M_{2}/m$ have the same order of magnitude $O(n^{2}p^{2})$. This led
to the following conjecture posed in~\cite{HaVu07}:
\begin{con}
\label{conjecture}
For all simple connected graphs~$G$:
$$
\frac {M_1}{n} \leqslant \frac{M_2}{m},
$$
and the bound is tight for complete graphs.
\end{con}
It was shown in~\cite{HaVu07} that this conjecture is not true in general
by finding a disconnected counterexample consisting of a six-vertex star and a triangle,
and a connected counterexample on $46$~vertices and $110$~edges.
Nevertheless, it was proven in~\cite{HaVu07} that
the conjecture holds for chemical graphs.
Further, it was proven in~\cite{VuGr07}
that the conjecture holds for trees (with equality attained for and only for stars),
while in~\cite{Li08} it was proven that the conjecture holds for connected unicyclic graphs
(with equality attained for and only for cycles).\rz

Our goal here is twofold, both in favor of a conjecture and against
it:
\begin{enumerate}[($i$)]
\item We show that the expressions $M_1/n$ and $M_2/m$ are both
bounded with $\frac{4m^2}{n^2}$ from below and with $\frac{\Delta
M_1}{2 m}$ from above, with equality attained for and only for
regular graphs. We also establish lower bounds for variable Zagreb
indices.

\item We show that for any fixed number $k\geq 2$,
there exists a connected graph with $k$~cycles for which
$M_1/n>M_2/m$ holds, effectively showing that the conjecture cannot
hold unless there exists some kind of limitation on the number of
cycles or the maximum vertex degree in a graph.  In particular, we
prove that the conjecture holds for subdivision graphs.
\end{enumerate}

\section{Common lower and upper bounds}

The following two theorems give sharp lower bounds for $M_1$
and~$M_2$. Recall that for a graph with $n$ vertices and $m$ edges,
the average value of vertex degrees is $2m / n$.

\begin{thm}
\label{first} It holds that $M_1 \geqslant \frac{4 m^2}{n}$.
The equality is attained if and only if graph is regular.
\end{thm}

\begin{proof}
We use the Cauchy-Schwartz inequality on vectors
$(d_1, d_2, \ldots, d_n)$ and $(1, 1, \ldots, 1)$ to get
$$
M_1 \cdot n = \left ( d_1^2 + d_2^2 + \ldots + d_n^2 \right ) \left
( 1^2 + 1^2 + \ldots + 1^2 \right ) \geqslant (d_1 \cdot 1 + d_2
\cdot 1 + \ldots + d_n \cdot 1)^2 = (2m)^2.
$$
Equality holds if and only if $d_1 = d_2 = \ldots = d_n$, namely
if and only if $G$ is regular.
\end{proof}

\begin{lemma}
\label{le-xlnx}
For positive real numbers $x_1, x_2, \ldots, x_n$ the following
inequality holds:
\begin{equation}
\label{eq-ineq}
x_1 \ln x_1 + x_2 \ln x_2 + \ldots + x_n \ln x_n \geqslant \left (
x_1 + x_2 + \ldots + x_n \right ) \ln \frac{x_1 + x_2 + \ldots +
x_n}{n}.
\end{equation}
\end{lemma}

\begin{proof}
The function $f (x) = x \ln x$ is strictly convex on interval $(0,+\infty)$,
since its second derivative $f'' (x) = \frac{1}{x}$ is positive.
The inequality~(\ref{eq-ineq}) follows directly from the Jensen's inequality~\cite{HaLP98}
$$
\frac {f (x_1) + f (x_2) + \ldots + f (x_n)}{n} \geqslant f \left
(\frac {x_1 + x_2 + \ldots + x_n}{n} \right ).
$$
Equality holds in~(\ref{eq-ineq}) if and only if all $x_i$ are equal.
\end{proof}

\begin{thm}
\label{second} It holds that $M_2 \geqslant \frac{4 m^3}{n^2}$.
The equality is attained if and only if graph is regular.
\end{thm}

\begin{proof}
First we use the inequality between the arithmetic and the geometric mean:
$$
\frac{M_2}{m} = \frac{\sum_{(i, j) \in E} d_i d_j}{m} \geqslant
\sqrt[m]{ \prod_{(i,j)\in E}^{\phantom{n}} d_id_j} =
\sqrt[m]{ \prod_{i = 1}^n d_i^{d_i} }.
$$
Since $d_i \geqslant 1$, we take the natural logarithm of both
sides to get
$$
\ln \frac{M_2}{m} \geqslant \frac{1}{m} \sum_{i=1}^n d_i\ln d_i.
$$
Then from Lemma~\ref{le-xlnx} we get:
$$
\ln \frac{M_2}{m} \geqslant \frac{1}{m}\,
\left(\sum_{i=1}^n d_i\right) \ln\left(\frac{\sum_{i=1}^n d_i}n\right)
=\frac1m\,2m \ln \frac{2m}{n}
= 2\ln \frac{2m}{n},
$$
and finally
$$
M_2\geqslant \frac{4m^3}{n^2}.
$$
Equality holds if and only if $d_1 = d_2 = \ldots = d_n$,
i.e., if and only if $G$ is regular.
\end{proof}

From two previous theorems, we see that the expressions from Conjecture~\ref{conjecture}
have common sharp lower bound:
$$
\frac{4m^2}{n^2} \leqslant \frac{M_1}n \quad \mbox{ and } \quad
\frac{4m^2}{n^2} \leqslant \frac{M_2}m.
$$

Next we show that these expressions also have common sharp upper bound.

\begin{prop}
\label{upper} Let $\Delta$ be the maximum vertex degree in~$G$. Then
$$
\frac{M_1}{n} \leqslant \frac{\Delta M_1}{2 m} \quad \mbox{
and } \quad \frac {M_2}{m} \leqslant \frac{\Delta M_1}{2 m}.
$$
Equality is attained simultaneously in both inequalities
if and only if $G$ is regular.
\end{prop}

\begin{proof}
The first inequality is equivalent to the obvious inequality $2m \leqslant \Delta n$,
while the second inequality is equivalent to $2 M_2 \leqslant \Delta M_1$.
Then
$$
2M_2 = \sum_{(i,j)\in E} 2d_id_j \leqslant \sum_{(i,j)\in E} (d_i^2
+ d_j^2) = \sum_{i\in V} d_i\cdot d_i^2 \leqslant \sum_{i\in V}
\Delta d_i^2 = \Delta M_1.
$$
Equality is attained in both inequalities simultaneously if and only
if $d_i=\Delta$ for every $1 \leqslant i \leqslant n$, i.e., if and
only if $G$ is regular.
\end{proof}

Now, using the upper bound on $M_1$ from \cite{Das04}
(where $\Delta$ is the maximum, while $\delta$ is the minimum vertex degree):
$$
M_{1}\leq m\left(\frac{2m}{n-1} + \frac{n-2}{n-1}\Delta
                +(\Delta-\delta)\left(1-\frac{\Delta}{n-1}\right)\right),
$$
with equality if and only if $G$ is a star graph or a regular graph
or $K_{\Delta+1}\cup(n-\Delta-1)K_{1}$,
we see that the expressions $M_1/n$ and $M_2/m$
also have common upper bound in terms of $n$, $m$, $\Delta$ and $\delta$:
\begin{eqnarray*}
\frac{M_1}{n} &\leqslant&
\frac{\Delta}2\left(\frac{2m}{n-1} + \frac{n-2}{n-1}\Delta
                +(\Delta-\delta)\left(1-\frac{\Delta}{n-1}\right)\right), \\
\frac {M_2}{m} &\leqslant&
\frac{\Delta}2\left(\frac{2m}{n-1} + \frac{n-2}{n-1}\Delta
                +(\Delta-\delta)\left(1-\frac{\Delta}{n-1}\right)\right).
\end{eqnarray*}
Equality is attained simultaneously in above inequalities if and
only if $G$ is regular. \rz

These indices have been generalized to variable first and second
Zagreb indices defined as
$$
{}^\lambda M_1 = \sum_{i = 1}^n d_i^{\ 2 \lambda} \quad \mbox{ and }
\quad {}^\lambda M_2 = \sum_{(i, j) \in E}^n (d_i d_j)^\lambda
$$

More results about comparing variable Zagreb indices can be found in
\cite{Vu07} and \cite{VuGr08}. For $2 \lambda \geqslant 1$, we
define $p = 2 \lambda$ and $q = \frac{2 \lambda}{2 \lambda - 1}$ in
order to establish relation $\frac{1}{p} + \frac{1}{q} = 1$. Now we
use H\"{o}lder inequality~\cite{HaLP98} on vectors $(d_1, d_2,
\ldots, d_n)$ and $(1, 1, \ldots, 1)$ to get
$$
\left ( \sum_{i = 1}^n d_i^{\ p} \right )^{1 / p} \cdot \left (
\sum_{i = 1}^n 1^q \right )^{1 / q} \geqslant \sum_{i = 1}^n (d_i
\cdot 1).
$$
Next, raise each side of equation to the power of $2 \lambda$
$$
\left ( \sum_{i = 1}^n d_i^{\ 2 \lambda} \right ) \cdot n^{2\lambda
- 1} \geqslant (2m)^{2 \lambda}.
$$
The last inequality is equivalent with
$$
{}^\lambda M_1 \geqslant n \left (\frac {2m}{n} \right )^{2
\lambda}.
$$

For the variable second Zagreb index and every $\lambda \geqslant 0$
it holds
$$
\frac{{}^\lambda M_2}{m} = \frac{\sum_{(i, j) \in E} (d_i
d_j)^{\lambda}}{m} \geqslant \sqrt[m]{ \prod_{(i,j)\in
E}^{\phantom{n}} (d_i d_j)^\lambda} = \sqrt[m]{ \prod_{i = 1}^n
d_i^{\lambda d_i} }.
$$
We can use the same technique as in the proof of Theorem
\ref{second} and get lower bound:
$$
{}^\lambda M_2 \geqslant m \left (\frac {2m}{n} \right )^{2
\lambda}.
$$

Also, we have similar upper bounds for variable Zagreb indices:
$$
\frac{{}^\lambda M_1}{n} \leqslant \frac{\Delta \cdot {}^\lambda
M_1}{2 m} \quad \mbox{ and } \quad \frac {{}^\lambda M_2}{m}
\leqslant \frac{\Delta \cdot {}^\lambda M_1}{2 m}.
$$

\section{Counterexamples}

Let $C (a, b)$ be a graph that is composed of $(a + 1)$-vertex star
with exactly $b$ triangles attached in line at arbitrary leaf (see Figure~1).
If triangles have vertex labels $v_i, u_i, w_i$, where $1
\leqslant i \leqslant k$, then there exist edges $u_i v_{i + 1}$ for
every $1 \leqslant i \leqslant k - 1$, and vertex $v_1$ is connected
with an arbitrary leaf of star $S_{a + 1}$. \rz
\begin{figure}[ht]
  \center
  \includegraphics [width = 10cm]{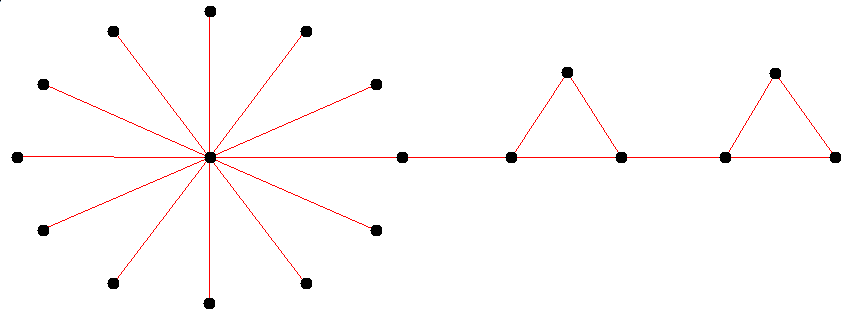}
  \caption { \textit{ The bicyclic counterexample $C (12, 2)$ with $19$ vertices }}
\end{figure}

Assume that $a \geqslant 3$ and $b \geqslant 1$. It is obvious that
the number of vertices of $C (a, b)$ is $n = a + 3b + 1$ and the
number of edges is $m = a + 4b$. Also note that $C (a, b)$ has
exactly $b$ cycles. \rz

In $C (a, b)$ there is one vertex of degree $a$ and $a - 1$ pendent
vertices. Every triangle has vertex degrees $3, 3, 2$, except for
the last one which has $3, 2, 2$. Now we can calculate the first
Zagreb index:
$$
M_1 (C (a, b)) =
a^2 + (a - 1) \cdot 1^2 + 2^2 + (b - 1) (3^2 + 3^2 + 2^2) + (3^2 + 2^2 + 2^2)
= a^2 + a + 22 b - 2.
$$

The weight of $a - 1$ pendent edges is equal to $a \cdot 1$, while
every triangle has weights $9, 6, 6$, except for the last one which
has $6, 6, 4$. The edges connecting triangles have weight $9$, and
therefore,
$$
M_2 (C (a, b)) = a \cdot 1 \cdot (a - 1) + 2 \cdot a + 2 \cdot 3 +
(9 + 6 + 6)(b - 1) + (6 + 6 + 2) + 9(b - 1)
= a^2 + a + 30 b - 8.
$$

The Conjecture \ref{conjecture} is equivalent to
$M_2 \cdot n - M_1 \cdot m \geqslant 0$,
which for the graph $C (a, b)$ yields:
$$
(a^2 + a + 30b - 8) (a + 3b + 1) - (a^2 + a + 22b - 2)(a + 4b) \geqslant 0.
$$
i.e.,
\begin{equation}
\label{eq-ab}
a^2 (1 - b) + a (7 b - 5) + (2 b^2 + 14 b - 8) \geqslant 0.
\end{equation}

Next, fix the number of cycles $b \geqslant 2$.
The left-hand side of~(\ref{eq-ab}) is a quadratic function in~$a$.
Since the coefficient of $a^2$ is negative, and the discriminant
$$
D = (7 b - 5)^2 - 4 (1 - b)(2 b^2 + 14 b - 8) = 8 b^3 + 97 b^2 - 158 b + 57
$$
is greater than zero for $b \geqslant 2$, we get that
the left-hand side value of~(\ref{eq-ab}) is negative for
\begin{equation}
\label{eq-alower}
a > \frac {-(7b - 5) + \sqrt{D}}{2 (1 - b)}.
\end{equation}

Thus, each value of~$a$ satisfying~(\ref{eq-alower}) yields
a counterexample to Conjecture~\ref{conjecture} with $b$~cycles.
In particular, for $b=2$ we get that any $a\geqslant 12$
yields a counterexample to the conjecture and
the smallest counterexample of this form is shown in Figure~1.

\section{Conclusion}

From the previous section it is evident that the
Conjecture~\ref{conjecture} cannot hold unless there exists some
kind of limitation on either the maximum vertex degree or the number
of cycles in a graph. This {\em limitation} may be implicitly given,
as it becomes evident from the following example. \rz

The subdivision graph $S (G)$ of a graph $G$ is obtained by
inserting a new vertex of degree two on each edge of $G$. If G has
$n$ vertices and $m$ edges, then $S (G)$ has $n + m$ vertices and $2
m$ edges. Clearly, $S (G)$ is bipartite.
\begin{thm}
Let $S (G)$ be a subdivision graph of $G$. Then,
$$
\frac {M_1 (S (G))}{n + m} \leqslant \frac {M_2 (S (G))}{2 m},
$$
with equality if and only if $G$ is a regular graph.
\end{thm}

\begin{proof}
The vertex degrees of~$G$ remain the same in the subdivision graph $S(G)$,
while the new vertices have degree two. Thus,
$$
M_1 (S (G)) = M_1 (G) + 2^2 \cdot m.
$$
Every edge $(i, j)$ of $G$ is subdivided in two parts with weights
$2 d_i$ and $2 d_j$. Therefore,
$$
M_2 (S (G)) = \sum_{(i, j) \in E} (2 d_i + 2 d_j) =
2 \sum_{i = 1}^n d_i^2 = 2M_{1}(G).
$$
Using these formulas, we get that
the inequality
$\frac{M_{1}(S(G))}{n+m}\leqslant \frac{M_{2}(S(G))}{2m}$
is equivalent to
$\frac {M_{1}(G) + 4m}{n + m} \geqslant \frac {2 M_{1}(G)}{2m}$,
i.e., to $M_{1}(G)\geqslant \frac{4m^{2}}n$,
which is true by Theorem~\ref{first}.
The case of equality easily follows.
\end{proof}

\end{document}